\documentclass[12pt]{article}
\usepackage{amsfonts}
\usepackage{graphicx}
\usepackage{subfig}
\usepackage{latexsym}
\usepackage{tabularx,array}
\usepackage{url}
\usepackage{amsmath,amsthm, dsfont, bm,bbm,xcolor}
\usepackage{mathtools}
\mathtoolsset{showonlyrefs}

\newcommand{\R}{\mathbb{R}}

\newcommand{\E}{\mathbb{E}}
\newcommand{\N}{\mathbb{N}}

\renewcommand{\P}{\mathbb P}

\newcommand{\al}{\alpha}

\newcommand{\ep}{\varepsilon}

\newcommand{\Indi}[1]{\mathbbm{1}_{#1}}

\newcommand{\rv}{\color{black}}

\newcommand{\bee}{\begin{equation}}
\newcommand{\eee}{\end{equation}}
\newcommand{\beea}{\begin{array}}
\newcommand{\eeea}{\end{array}}

\renewcommand{\theequation}{\arabic{section}.\arabic{equation}}

\theoremstyle{plain}
\newtheorem{prop}{Proposition}[section]

\newtheorem{theo}[prop]{Theorem}

\newtheorem{lemma}[prop]{Lemma}

\theoremstyle{definition}

\newtheorem{rem}[prop]{Remark}

\begin{document}

\title{Optimal estimation of local time and occupation time measure for an $\al$-stable L\'evy process 
\thanks{The authors
gratefully acknowledge financial support of ERC Consolidator Grant 815703
``STAMFORD: Statistical Methods for High Dimensional Diffusions''.}} 
\author{Chiara Amorino\thanks{Department
of Mathematics, University of Luxembourg, 
E-mail: chiara.amorino@uni.lu.}  \and
Arturo  Jaramillo\thanks{Center of Research in Mathematics (CIMAT), Guanajuato, 
E-mail: jagil@cimat.mx.} \and
Mark Podolskij\thanks{Department
of Mathematics, University of Luxembourg,
E-mail: mark.podolskij@uni.lu.}}

\maketitle

\begin{abstract}
\noindent 
{We} present a novel theoretical result on estimation
of local time and occupation time measure of
{\rv an $\al$-stable} L\'evy process with $\alpha \in (1,2)$. Our approach is based upon computing the 
conditional expectation of the desired quantities given high frequency data, which is an $L^2$-optimal statistic by construction. We prove the corresponding stable central limit theorems and discuss a statistical application. In particular, this work extends the results of \cite{IP21}, which investigated the case of the Brownian motion.

\ \

\bigskip


\end{abstract}

\section{Introduction} \label{sec1}
\setcounter{equation}{0}
\renewcommand{\theequation}{\thesection.\arabic{equation}}

\subsection{The setting and overview} 
In this paper we consider a {\rv pure jump} $\alpha$-stable L\'evy process 
{$X=\{X_t\}_{t\geq 0}$} with $\alpha \in (1,2)$, defined on a filtered probability space $(\Omega,\mathcal{F},{\{\mathcal{F}_{t}\}_{t\geq 0}},\mathbb{P})$. The law of the process $X$ is uniquely determined by the  characteristic function
\begin{align}\label{eq:characteristifcuntion}
\E[\exp(\textbf{i}\xi X_{t})]
  &=\exp({\rv -\bm{i}\xi t\eta -\sigma t |\xi|^{\alpha}(1-\bm{i}\beta\tan(\frac{\pi\alpha}{2})\text{sgn}(\xi))})
\end{align}
{\rv for $\xi\in\R$, $t\geq0$ $\eta\in\R$, $\beta\in[-1,1]$ and $\sigma>0$.} We will focus on estimation of occupation time measure and the local time of $X$ given high frequency data {$\{X_{i/n}\}_{1\leq i \leq \lfloor nT \rfloor}$, with $T>0$} being fixed and $n\to \infty$. We recall that for  
$a<b$ the occupation time of $X$ at the set $(x, \infty)$
over the interval $[a,b]$, denoted by $O_{[a,b]}(x)$, is defined as 
$$O_{[a,b]}(x):= \int_a^b \Indi{(x, \infty)}(X_s) ds.$$
The local time of $X$ at point $x\in\R$ over the interval $[a,b]$ denoted as 
$L_{[a,b]}(x)$ is defined implicitly via the
\textit{occupation density formula}:
\begin{align}\label{eq:Oabtolocaltime}
O_{[a,b]}(x)= \int_x^{\infty} L_{[a,b]}(y) dy \qquad \mathbb{P}-a.s.
\end{align}
(if it exists). Throughout the paper we use the abbreviation $L_t:=L_{[0,t]}$ and $O_t:=O_{[0,t]}$.  The existence and smoothness properties of local times of stochastic processes have been extensively studied in the 70's and 80's; we refer to articles \cite{B85,B88,B69,B70,GH80,GH72} among many others. In particular, in the setting of pure jump $\al$-stable L\'evy processes, the local time exists only for $\al \in (1,2)$ (cf. \cite{K69}). Furthermore, {\rv there exists a version of local time, which is continuous 
in space and time.} Indeed, according to  \cite[Theorems 2 and 3]{B88} and
\cite[Theorem 4.3]{hold}, {\rv there exists a version of the local time with } the following properties:  
\bee \label{smoothn}
L_t(\cdot) \text{ is $\P$-a.s. locally H\"older continuous with index } (\al-1)/2-\ep,
\eee
for any $\ep \in (0, (\al-1)/2)$. Moreover, for all $\varepsilon\in(0,1-1/\alpha)$ and $T>0$, there exists a deterministic constant {\rv $C_T>0$}, such that 
\begin{align}\label{eq: 1.4.5}
\sup_{0\leq s\leq t\leq T}|L_{t}(x)-L_{s}(x)|\leq {\rv C_T}|t-s|^{(1-1/\alpha)-\varepsilon}\ \ \ \ \ \   \mathbb{P}\text{-a.s.}
\end{align}
In particular, $L_{(\cdot)}(x)$  is $\P$-a.s. locally H\"older continuous with index $(1-1/\alpha)-\varepsilon$ for any $\ep \in (0, 1- 1/\alpha)$. {\rv Throughout the paper we consider the aforementioned H\"older continuous version of the local time.}

The goal of this article is to introduce the $L^2$-optimal estimators of $L_{t}(x)$ and $O_{t}(x)$, and study their asymptotic properties. For this purpose we define
the $\sigma$-algebra 
$\mathcal{A}_n:=\sigma(\{X_{i/n}\ ;\ i\in\N\})$
and construct the estimators via
\begin{align}
S_{L,t}^{n}{\rv (x)}
  =\E[L_{t}(x)\ |\mathcal{A}_n]
\ \ \ \text{and} \ \ \
S_{O,t}^{n}{\rv (x)}
  =\E[O_{t}(x)\ |\mathcal{A}_n]\label{eq:testim}.
\end{align}
By construction $S_{L,t}^{n}$ (resp. $S_{O,t}^{n}$) is an $L^2$-optimal approximation of 
$L_t(x)$ (resp. $O_t(x)$). We will show a functional \textit{stable convergence} {\rv for both statistics, which exhibit a mixed normal limit}. Our main tool is Jacod's stable central limit theorem for partial sums of high frequency data stated in \cite{Jacod}.

\subsection{Related literature} \label{sec1.2}
Functional limit theorems for estimators of local times have been studied for numerous stochastic models. In the setting of $\al$-stable L\'evy processes {and} related models, such estimators often take the form $G(x,\phi)^n=\{G(x,\phi)_t^n;\ t\geq 0\}$, with
\begin{align} \label{statistic}
G(x,\phi)_t^n=n^{\frac{1}{\alpha}-1}\sum_{ i = 1}^{\lfloor nt\rfloor}\phi\left(n^{\frac{1}{\alpha}}(X_{\frac{i-1}{n}}-x)\right),
\end{align}
where $\phi\in L^{1}(\R)$. Consistency and asymptotic mixed normality for statistics $G(x,\phi)^n$ and related functionals have been investigated in  \cite{B86,J98} in the setting of Brownian motion and continuous stochastic differential equations. Some optimality results for estimation of the local time and occupation time measure of a Brownian motion can be found in \cite{IP21}. 
Estimation errors for occupation time functionals of stationary Markov processes have been studied in \cite{AC18}. Limit theorems for statistics of the form \eqref{statistic} in the case of the fractional Brownian motion are discussed in \cite{JNP21,J04,J06,PR18}, although the complete weak limit theory is far from being understood. 

{\rv While estimation of local time and occupation time measure has an interest in its own right, accurate estimation of these objects can be useful for related statistical problems. For example, non-parametric estimators of the diffusion coefficient in a continuous stochastic differential equation often involve local times in the mixed normal limit, see e.g. \cite{F93}. In a similar spirit, local times and mixed normal local times appear as fundamental limits for additive functionals of a variety of Gaussian processes (see Papanicolaou, Stroock and Varadhan \cite{PaStVar},  Jaramillo, Nourdin, Nualart and Peccati \cite{JNP212} as well as Minhao Hong, Heguang Liu and Fangjun Xu \cite{HoLiuXu}), thus serving as simplifying probabilistic models for otherwise complex probabilistic objects. Efficient estimation of local times and occupation times is very beneficial for statistical inference in this framework.}

Our main result about the asymptotic theory for local times is mostly related to the articles \cite{J04,J06,J08,IP21, R91}. The paper  \cite{J04} {addresses consistency in the case of} the linear fractional stable motion, a class which in particular includes stable L\'evy processes. 

\begin{theo}[Theorem 4 in \cite{J04}]\label{thm:Jeganathan}
Suppose that $\phi:\R\rightarrow \R$ is a function satisfying $\phi,\phi^2\in L^1(\R;\R)$. Then, for every $t>0$, 
\begin{align*}
G(x,\phi)_t^n
  &\stackrel{L^{2}(\Omega)}{\rightarrow} \int_{\R}\phi(y)dy\cdot L_{t}(x) \qquad
  {\rv \text{as } n\to \infty}.
\end{align*}
\end{theo}

\noindent
In \cite{J06,J08,R91} the authors prove the asymptotic mixed normality for continuous version of the functional $ G(x, \phi)^n$, but only in the \textit{zero energy} setting, i.e. $\int_{\R}\phi(y)dy=0$.  
From statistical point of view,  this case is a less interesting one as we would like to use statistics of the type  \eqref{statistic} as an estimator of the local time $L_t(x)$. More importantly, in the setting $\int_{\R}\phi(y)dy \not=0$ the methods {\rv developed} in  \cite{J06,J08} do not apply.

{\rv We will  use  stable convergence theorems for high frequency statistics {\rv introduced} in \cite{Jacod} to show asymptotic mixed normality of standardised versions of $S_{L,t}^n(x)$
and $S_{O,t}^n(x)$. These results are strongly related to an earlier work \cite{IP21}, which considers the same problem in the setting of the Brownian motion. While their limit theorems are also based on the general results of \cite{Jacod}, the technical aspects of the proof are more involved in the case of pure jump L\'evy processes. The details will be highlighted in the proof section.}

\subsection{Outline of the paper} \label{sec1.3}
The rest of the paper is organised as follows. Section \ref{sec2} presents some preliminaries, main results and an application. Proofs of main results are collected in Section \ref{sec3}. Some technical statements are proved in Section \ref{sec4}.

\section{Preliminaries and main results} \label{sec2}
In this section we present the spectral representation of  local times, which will be useful in the sequel. Furthermore, we introduce the notion of stable convergence and establish the asymptotic theory for the estimators $S^n_{L}$ and $S^n_{O}$.
\subsection{Local times and stable convergence}\label{sec:localtime}

The analysis of occupation times and local times {is } an integral part of the {theory of  stochastic processes}, which found manifold applications in probability during the past decades. We recall that, for $t> 0$ and $x \in \R$, the local time of {\rv the $\alpha$-stable L\'evy process} $X$, 
$\alpha\in (1,2)$, up to time $t$ at $x$ can be formally defined as $L_t(x) := \int_0^t \delta_0(X_s - \lambda) ds$, where $\delta_0$ denotes the Dirac delta function. A rigorous definition of the local time is obtained by replacing $\delta_0$ by the Gaussian kernel $\phi_\epsilon(x) := (2 \pi \epsilon)^{- \frac{1}{2}} \exp(- \frac{1}{2 \epsilon}x^2)$ and taking the limit in probability for $\epsilon \rightarrow 0$.
For our purposes we will systematically use the following spectral representation for $L_{t}(x)$. {\rv The proof of the representation gathered in the lemma below can be found in e.g.  \cite[Proposition 11]{J04}, for completeness we give the proof of the boundedness of the moments. It can be found in Section \ref{sec4}.} 

\begin{lemma}\label{l:lemmaFourierLT}
For every $t>0$ and $x\in\R$, the sequence
$$\left\{\int_{[-m,m]}\int_{0}^te^{-\textbf{i} \xi (X_{s}-x)} {\rv dsd\xi}  \right\}_{ m\in\N}$$
converges in $L^{2}(\Omega)$. The limit as $m$ tends to infinity, which we will denote by 
$$\int_{\R}\int_{0}^te^{-\textbf{i} \xi (X_{s}-x)}{\rv dsd\xi}$$
satisfies 
\begin{align}\label{eq:FourierrepLt}
L_{t}(x)
  &=\frac{1}{2\pi}\int_{\R}\int_{0}^te^{-\textbf{i} \xi (X_{s}-x)}{\rv dsd\xi}.
\end{align}
{\rv Moreover, for any $p \in \mathbb{N}$, $\E[(L_t(x))^p] < \infty$.}
\end{lemma}
\noindent

{\rv Recall furthermore that, by the definition of the occupation time $O_t(x)$, for any $p \in \mathbb{N}$  $\E[(O_t(x))^p] < t^p$.}

In what follows we will use the notion of \textit{stable convergence}, which is originally due to Renyi \cite{R63}. Let $(\mathcal{S},\delta)$ be a Polish space. Let  {$\{Y_n\}_{n\geq 1}$} be a sequence of $\mathcal{S}$-valued and $\mathcal{F}$-measurable random variables defined on $(\Omega,\mathcal{F},\mathbb{P})$ and $Y$ a random variable defined on an enlarged probability space $(\Omega',\mathcal{F}',\mathbb{P}')$. We say that $Y_n$ converges stably to $Y$, if and only if for any continuous bounded function $g:\mathcal{S} \to \R$ and any bounded
$\mathcal F$-measurable random variable $F$, it holds that
\bee \label{gstab}
\lim_{n \to \infty} \E[Fg(Y_n)] = \E'[Fg(Y)].
\eee
In this case we write $Y_n \stackrel{Stably}{\rightarrow} Y$. In this paper we deal with the space 
{\rv of c\'adl\'ag functions $D([0,T])$ equipped with the Skorokhod $J_1$-topology}.

\subsection{Main results}
The following proposition provides an explicit  expression of the statistics $S_{L,t}^{n}$ and $S_{O,t}^{n}$. {\rv The  proof 
of Proposition \ref{l: estim local}} is based upon the Markov property of $X$ and the linearity in time of our objects. 

\begin{prop}\label{l: estim local}
For $i\geq 1$, define the increments $\Delta_i^n X := X_{\frac{i}{n}} - X_{\frac{i - 1}{n}}$. Consider  the function $f:\R^{2}\rightarrow \R$ given by $f(x, y):= \E[L_{[0,1]}(x) | X_1 = y]$ and $F(x,y):=\int_{x}^{\infty}f(r,y)dr$. Then we obtain the identities
\begin{align}
S_{L,t}^n(x)
  &= n^{\frac{1}{\alpha} - 1} \sum_{i = 1}^{\lfloor nt \rfloor} f\left(n^{\frac{1}{\alpha}}(x - X_{\frac{i - 1}{n}}), n^{\frac{1}{\alpha}} \Delta_i^n X\right) + 
\mathcal{E}_{L,t}^n(x)\label{eq:SLexact} \\
S_{O,t}^n(x)
  &= \frac{1}{n} \sum_{i = 1}^{\lfloor nt \rfloor} F\left(n^{\frac{1}{\alpha}}(x - X_{\frac{i - 1}{n}} ), n^{\frac{1}{\alpha}} \Delta_i^n X\right) + 
\mathcal{E}_{O,t}^n(x)\label{eq:SOexact}
\end{align}
where $\mathcal{E}_{L,t}^n(x)$ and $\mathcal{E}_{O,t}^n(x)$ are defined as 
\begin{align*}
\mathcal{E}_{L,t}^n(x)
  &:=\E\left[L_{[\lfloor nt \rfloor/n, t]}(x) | \mathcal{A}_n\right] \\[1.5 ex]
  \mathcal{E}_{O,t}^n(x)
  &:= \E\left[O_{[\lfloor nt \rfloor/n, t]}(x) | \mathcal{A}_n\right].
\end{align*}
Moreover, the processes $\{n^{\frac{1}{2}(1-\frac{1}{\alpha})}\mathcal{E}_{L,t}^n(x)\}_{t\geq 0}$ and $\{n^{\frac{1}{2}(1+\frac{1}{\alpha})}\mathcal{E}_{O,t}^n(x)\}_{t\geq 0}$ converge to zero in probability
uniformly on compact intervals {\rv as  $n\to \infty$}. 
\end{prop}

The proof of Proposition \ref{l: estim local} will be given in Section \ref{sec4}.
The next theorem is a functional limit result for the error of the approximation of the occupation and local time by their $L^{2}$-optimal estimators. 

\begin{theo}\label{th: main local}
Fix $x\in \R$ and define the processes 
\begin{align*}
W_{L,t}^n
  &:= n^{\frac{1}{2}(1 - \frac{1}{\alpha})}(S_{L,t}^n(x) - L_{t}(x)) \\
W_{O,t}^n
  &:= n^{\frac{1}{2}(1 + \frac{1}{\alpha})}(S_{O,t}^n(x) - O_{t}(x)).
\end{align*}
Then we obtain the functional stable convergence with respect to {\rv $J_1$} topology:
\begin{align*}
W_{L}^n
  \stackrel{Stably}{\rightarrow} k_L B_{L(x)}, \qquad W_{O}^n
  \stackrel{Stably}{\rightarrow} k_O B_{L(x)}, \qquad
  {\rv \text{as } n\to \infty,} \\
\end{align*}
where $B$ is a standard Brownian motion defined on an extended space and independent of $\mathcal{F}$. The constants $k_L$ and $k_O$ are defined as
\begin{align*}
k_L^2
  &:= \int_{\R} \E[(\E[L_{1}(y) | X_1] - L_{1}(y))^2] dy\\
k_O^2
  &:= \int_{\R} \E \big[(\E[O_{1}(z) | X_1 ] - O_{1}(z))^2 \big] dz.
\end{align*}
\end{theo}
\begin{rem}
By the Dambis-Dubins-Schwarz theorem {\rv(cf. \cite[Theorem 1.6, Section 5.1]{RY99})},  {\rv Theorem \ref{th: main local} implies 
the stable convergence
\begin{align*}
  W_{L}^n
  \stackrel{Stably}{\rightarrow} k_L {\int_{0}^{\cdot}\sqrt{L_{s}(x)}B(ds)}, \qquad
  W_{O}^n
  \stackrel{Stably}{\rightarrow} k_O {\int_{0}^{\cdot}\sqrt{L_{s}(x)}B(ds)}.
\end{align*}
as  $n\to \infty$.
}
\qed
\end{rem}
\begin{rem}
The statement above can be directly used to construct confidence regions for $L_{t}(x), O_{t}(x) $ if the law of the L\'evy process $X$ is known. 
Indeed, by properties of stable convergence, it holds for any fixed
$t>0$ {\rv as  $n\to \infty$}:
\bee
\frac{n^{\frac{1}{2}(1 - \frac{1}{\alpha})} \left(  S_{L,t}^n(x) -   L_{t}(x)  \right)}{ k_L \sqrt{S_{L,t}^n(x)}} \stackrel{d}{\rightarrow} \mathcal N(0,1)
\eee
and
\bee
\frac{n^{\frac{1}{2}(1 + \frac{1}{\alpha})} \left(  S_{O,t}^n(x) -   O_{t}(x)  \right)}{ k_O \sqrt{S_{{\rv L},t}^n(x)}} \stackrel{d}{\rightarrow} \mathcal N(0,1).
\eee

Asymptotic confidence sets for $L_{t}(x)$ and $O_{t}(x)$ readily follow from the above central limit theorem. \qed
\end{rem}

{\rv
\begin{rem}
The setting of {\rv $\alpha$-stable L\'evy processes,} $\alpha \in (1,2)$, is rather convenient as we often use self-similarity of $X$ in our arguments. However, it might not be a necessary assumption. We conjecture that similar results can be shown for \textit{locally $\alpha$-stable}  L\'evy processes although some bias {\rv effects} may appear in the limit theory. 
\qed 
\end{rem}
}

\noindent
From statistical point of view Theorem \ref{th: main local} provides lower bounds for estimation of the path functionals $L_{t}(x) $ and $O_{t}(x)$. In particular, it shows that statistics  $G(x,\phi)^n$ introduced in \eqref{statistic} do not produce $L^2$-optimal estimates. On the other hand, $G(x,\phi)^n$ can be computed from data even when the exact law of $X$ is unknown (use e.g. $\phi=\Indi{[-1,1]}$) in contrast to $L^2$-optimal statistics (recall that the functions $f$ and $F$ depend on the parameters of the stable distribution). We remark however that the limit theory for statistics  $G(x,\phi)^n$ is expected to be much more involved; see the proofs \cite{J98} for more details in the case of Brownian motion. Hence, we postpone the discussion to future research.

\section{Proof of main results}\label{sec3}
This section is devoted to the proof of the central limit theorem stated in the previous section. Throughout the proofs we denote by $C>0$ a generic constant, which may change from line to line.  In the proof of our main results the following lemma will be repeatedly used. Its proof can be found in Section \ref{sec4}.

\begin{lemma}
{\rv Let us introduce $\varphi_1(y) := \E[L^p_{[0,1]}(y)]$ and $\varphi_2(y) := \E[O^p_{[0,1]}(y)]$. Then, for} any $x \in \mathbb{R}$, $p \in \mathbb{N}$ and $i \in \left \{ 1, \ldots ,n \right \}$, the following identities hold true. 
\begin{enumerate}
    \item[(a)] $\E[L^p_{[\frac{i-1}{n}, \frac{i}{n}]}(x) | \mathcal{F}_{\frac{i - 1}{n}}] = n^{p (\frac{1}{\alpha} - 1)} {\rv \varphi_1 (n^{\frac{1}{\alpha}}(x - X_{\frac{i-1}{n}})) } $,
    \item[(b)] $\E[O^p_{[\frac{i-1}{n}, \frac{i}{n}]}(x) | \mathcal{F}_{\frac{i - 1}{n}}] = n^{-p} \, {\rv \varphi_2 (n^{\frac{1}{\alpha}}(x - X_{\frac{i-1}{n}})) } $.
\end{enumerate}
Moreover, define the functions {\rv $\varphi_3(y) := \E[ \E[L_{[0,1]}(y)| X_1]L_{[0,1]}(y)]$ and $\varphi_4(y) := \E[ \E[O_{[0,1]}(y)| X_1]O_{[0,1]}(y)]$ and recall that $f(x,y)=\E[L_{1}(x) | X_1 = y]$ and $F(x,y) = \int_x^\infty f(r,y) dr$.}
Then, 
\begin{enumerate}
    \item[(c)] $\E \left[f\left(n^{\frac{1}{\alpha}}(x - X_{\frac{i - 1}{n}}), n^{\frac{1}{\alpha}} \Delta_i^n X \right)
    L_{[\frac{i-1}{n}, \frac{i}{n}]}(x) | \mathcal{F}_{\frac{i - 1}{n}} \right]$
    $= n^{\frac{1}{\alpha} - 1} {\rv \varphi_3 (n^{\frac{1}{\alpha}}(x - X_{\frac{i-1}{n}})) } $,
    \item[(d)] $\E \left[F\left(n^{\frac{1}{\alpha}}(x - X_{\frac{i - 1}{n}}), n^{\frac{1}{\alpha}} \Delta_i^n X \right)
    O_{[\frac{i-1}{n}, \frac{i}{n}]}(x) | \mathcal{F}_{\frac{i - 1}{n}} \right]$
    $= n^{ - 1} {\rv \varphi_4 (n^{\frac{1}{\alpha}}(x - X_{\frac{i-1}{n}})) }$.
\end{enumerate}
\label{l: moments L and O}
\end{lemma}

{\rv Regarding the well-posedness of the elements mentioned above, it is worth noting that $\varphi_1$ and $\varphi_2$ are finite, as shown in Lemma \ref{l:lemmaFourierLT} and the comment below it. It is also easy to verify that $\varphi_3$ and $\varphi_4$ are finite. In particular, demonstrating the finiteness of $\varphi_3$ requires the use of Cauchy-Schwarz and Jensen's inequalities, in combination with the boundedness of the moments of $L_{1}$ as established in Lemma \ref{l:lemmaFourierLT}.}

\subsection{Proof of Theorem \ref{th: main local}}
\begin{proof}
Our argument is based on a martingale approach. We first deal with the estimation of the local time. Recall that $f(x,y)=\E[L_{1}(x) | X_1 = y]$. Because of Proposition \ref{l: estim local}
{\rv and the linearity of the local time (in time)}, we can write
\begin{equation}
W_{L,t}^n = n^{\frac{1}{2}(1 - \frac{1}{\alpha})}\left(S_{L,t}^n(x) - L_{[0,t]}(x)\right)
= \sum_{i = 1}^{\lfloor nt \rfloor} Z_{in}^{L}  - n^{\frac{1}{2}(1 - \frac{1}{\alpha})}N_{L,t}^{ n} + n^{\frac{1}{2}(1 - \frac{1}{\alpha})}\mathcal{E}_{L,t}^{n}(x),
\label{eq: Wtn start}
\end{equation}
where 
$$Z_{in}^L:= n^{\frac{1}{2}(\frac{1}{\alpha}-1)} f\left(n^{\frac{1}{\alpha}}(x - X_{\frac{i - 1}{n}}), n^{\frac{1}{\alpha}} \Delta_i^n X \right) - n^{\frac{1}{2}(1- \frac{1}{\alpha})} L_{[\frac{i-1}{n}, \frac{i}{n}]}(x)$$ 
and 
$$N_{L,t}^{ n}= L_{[\lfloor nt \rfloor /n, t]}(x).$$ 
The term $N_L^{ n}$ appears due to the edge effect and it does not contribute to the limit as stated in the following lemma, whose proof can be found in the next section {\rv and immediately follows from \eqref{eq: 1.4.5}.} 
\begin{lemma}
The process $N_L^{n}=\{N_{L,t}^{ n}\ ;\ t\geq0\}$ satisfies the following convergence  in probability uniformly over compact sets
$$\lim_n n^{\frac{1}{2}(1 - \frac{1}{\alpha})} N_L^{ n}  = 0. $$
\label{l: edge local}
\end{lemma}
\noindent Due to Proposition \ref{l: estim local} and Lemma \ref{l: edge local} we can write  \eqref{eq: Wtn start}  as 
$$W_{L,t}^n = \sum_{i = 1}^{\lfloor nt \rfloor} Z_{in}^L + o_{\mathbb{P}}(1).$$
We introduce here the notation $\E_{\frac{i-1}{n}}[\cdot ]$ for $\E[\cdot|  \mathcal{F}_{\frac{i-1}{n}}]$, which will be useful in the sequel.\\

\noindent  In the next step we will apply Theorem 3-2 of \cite{Jacod}. In our setting,
  $\E_{\frac{i-1}{n}}[Z_{in}^L] = 0$  because of {\rv Equation \eqref{eq: local intervals} below.} Then, it suffices to show
the following conditions:
\begin{equation}
\sum_{i = 1}^{\lfloor nt \rfloor} \E_{\frac{i-1}{n}}[|Z_{in}^L|^2] \xrightarrow{\mathbb{P}} k_L^2 L_{t}(x) \qquad \forall t \in [0, 1],
\label{eq: cond1 local}
\end{equation}
\begin{equation}
\sum_{i = 1}^{\lfloor nt \rfloor} \E_{\frac{i-1}{n}}[Z_{in}^L \Delta_i^n M]  \xrightarrow{\mathbb{P}} 0 \qquad \forall t \in [0, 1],
\label{eq: 3.15 Jacod}
\end{equation}
\begin{equation}
\sum_{i = 1}^{\lfloor nt \rfloor} \E_{\frac{i-1}{n}}\left[|Z_{in}^L|^2 \Indi{\{|Z_{in}^L| > \epsilon\}} \right]  \xrightarrow{\mathbb{P}} 0 \qquad \forall \epsilon > 0,
\label{eq: 3.13 Jacod}
\end{equation}
where the condition \eqref{eq: 3.15 Jacod} should hold for all square integrable continuous martingales $M$. {\rv We emphasise that we have summarised the two conditions (3.12) and (3.14) from 
\cite[Theorem 3-2]{Jacod} into our constraint \eqref{eq: 3.15 Jacod}}. Indeed, conditions (3.12) and (3.14) in \cite[Theorem 3-2]{Jacod} are formulated with respect to some continuous martingale $M$ which has to be chosen by the user according to the problem at hand. The convergence in \eqref{eq: 3.15 Jacod} holding for all bounded continuous martingales implies them both. It entails, in particular, that the continuous process $G$ of \cite[Theorem 3-2]{Jacod} is in our case identically zero. \\

\noindent We start by showing condition \eqref{eq: cond1 local}. Due to the identity $f(x, y)= \E[L_{[0,1]}(x) | X_1 = y]$  we can {\rv check} that
\begin{align}{\label{eq: conv squared L}}
 \sum_{i = 1}^{\lfloor nt \rfloor} \E_{\frac{i-1}{n}}[|Z_{in}^L|^2] = n^{\frac{1}{\alpha}-1}\sum_{i = 1}^{\lfloor nt \rfloor} \psi\left(n^{\frac{1}{\alpha}}(x - X_{\frac{i - 1}{n}})\right),   
\end{align}
where $\psi$ is given by 
$$\psi(q):=\E\big[  (\E[L_{[0, 1]}(q) | X_1 ] - L_{[0, 1]}(q))^2 \big].$$
Indeed, from the definition of $Z_{in}^L$ we have 
\begin{align*}
\sum_{i = 1}^{\lfloor nt \rfloor} \E_{\frac{i-1}{n}}[|Z_{in}^L|^2] = & \sum_{i = 1}^{\lfloor nt \rfloor}  \big[n^{\frac{1}{\alpha}-1} \E_{\frac{i-1}{n}}[f^2(n^{\frac{1}{\alpha}}(x - X_{\frac{i - 1}{n}}), n^{\frac{1}{\alpha}} \Delta_i^n X)] + n^{1 - \frac{1}{\alpha}} \E_{\frac{i-1}{n}}[L^2_{[\frac{i-1}{n}, \frac{i}{n}]}(x)] \\
& -2  \E_{\frac{i-1}{n}}[f(n^{\frac{1}{\alpha}}(x - X_{\frac{i - 1}{n}}), n^{\frac{1}{\alpha}} \Delta_i^n X)L_{[\frac{i-1}{n}, \frac{i}{n}]}(x)]  \big]\\
=:& I_1 + I_2 + I_3
\end{align*}
We start studying $I_1$. It is easy to see that
\begin{align*}
I_1 &= n^{\frac{1}{\alpha}-1} \sum_{i = 1}^{\lfloor nt \rfloor} \E_{\frac{i-1}{n}}[f^2(n^{\frac{1}{\alpha}}(x - X_{\frac{i - 1}{n}}), n^{\frac{1}{\alpha}} \Delta_i^n X)] \\
& = n^{\frac{1}{\alpha}-1} \sum_{i = 1}^{\lfloor nt \rfloor} g_1(n^{\frac{1}{\alpha}}(x - X_{\frac{i - 1}{n}})),
\end{align*}
with {\rv $g_1(q) := \E[f^2(q, X_1)]$}. From Lemma \ref{l: moments L and O}(a) the following is straightforward: 
\begin{align*}
I_2 & = n^{1 - \frac{1}{\alpha}} \sum_{i = 1}^{\lfloor nt \rfloor} \E_{\frac{i-1}{n}}[L^2_{[\frac{i-1}{n}, \frac{i}{n}]}(x)]  \\
& = n^{ \frac{1}{\alpha}-1} \sum_{i = 1}^{\lfloor nt \rfloor} g_2(n^{\frac{1}{\alpha}}(x - X_{\frac{i - 1}{n}})),
\end{align*}
{\rv where $g_2(q):=\E[L_{[0,1]}(q)^2]$.}
We are left to study $I_3$. From Lemma \ref{l: moments L and O}(c) we directly obtain 
$$I_3 = n^{ \frac{1}{\alpha}-1} \sum_{i = 1}^{\lfloor nt \rfloor} g_3(n^{\frac{1}{\alpha}}(x - X_{\frac{i - 1}{n}})),$$
with $g_3(q) := -2 \E[L_{[0,1]}(q) \E[L_{[0,1]}(q) |X_1]]$. 
Putting everything together we get \eqref{eq: conv squared L} with $\psi(q)= g_1(q) + g_2(q) + g_3(q)$.

To apply Theorem \ref{thm:Jeganathan} we need to check that $\psi, \psi^2 \in L^1(\R)$. By Minkowski and Jensen inequalities, $\psi(q) \le 4 \E[L_{[0, 1]}^2(q)]$. {\rv We therefore aim at showing that the function $\E[L_{[0, 1]}^2(\cdot )]$ belongs to $L^1(\R)$.} Let $\tau_q$ denote the first passage time of $X$ over the level $q$. Conditioning over $\tau_q$,  using the additivity of the local time, we deduce the inequality
$$\E[L_{[0, 1]}^2(q)] \le \E[L_{[0, 1]}^2(0)] \mathbb{P}(\tau_q < 1).$$ 
Moreover, by the Fourier representation of the local time (as stated in Lemma \ref{l:lemmaFourierLT}), $\E[L_{[0, 1]}^2(0)]$ is bounded. On the other hand, 
\begin{equation}
\int_0^\infty \mathbb{P}[\tau_q < 1] dq = \int_0^\infty \mathbb{P}[\sup_{s \le 1} X_s > q] dq = \E[\sup_{s \le 1} X_s].
\label{eq: prob passage time}
\end{equation}
We recall that $\E[(\sup_{s \le 1} X_s)^p]$ is bounded for any $p \in (0, \alpha)$ {\rv (see Corollary II.1.6 and Theorem II.1.7 in \cite{RY99})}. As $\alpha \in (1,2)$, this implies the boundedness of $\E[\sup_{s \le 1} X_s]$. Hence, by symmetry, $\psi \in L^1(\R)$ and $k_L^2 < \infty$. Applying the same argument it is easy to show that also $\psi^2\in L^1(\R)$, thanks to \eqref{eq: prob passage time} and the boundedness of $\E[L_{[0, 1]}^4(q)]$.
By Theorem \ref{thm:Jeganathan} we conclude that 
\[
\sum_{i = 1}^{\lfloor nt \rfloor} \E_{\frac{i-1}{n}}[|Z_{in}^L|^2] \xrightarrow{\mathbb{P}} k_L^2 L_{t}(x),
\]
where $k_L^2=\int_{\R} \psi(y) dy < \infty$.\\ 

\noindent In the next step we show condition \eqref{eq: 3.15 Jacod}.  Let $M$ be any continuous square integrable martingale and denote by $\mu$ (resp. $\widetilde{\mu}$) the random measure (resp. compensated random measure) associated with
the pure jump L\'evy process $X$.  The martingale representation theorem for jump measures investigated in Lemma 3 (ii) and Theorem 6 of \cite{Repr} implies that $Z_{in}^L$ has an integral representation
\[
Z_{in}^L = \int_{\frac{i-1}{n}}^{\frac{i}{n}} \int_{\R} \eta_i^n(x,t) 
\widetilde{\mu}(dx,dt) 
\]
for some predictable square integrable process $\eta_i^n$.  
Using that the covariation between any continuous martingale and any pure jump martingale is zero, we conclude that 
$$\E_{\frac{i-1}{n}}[Z_{in}^L \Delta_i^n M] = 0.$$
Thus, we obtain \eqref{eq: 3.15 Jacod}. \\

\noindent Finally, we show condition \eqref{eq: 3.13 Jacod}.  {\rv Cauchy-Schwarz} inequality ensures that
{\rv
$$\E_{\frac{i-1}{n}}[|Z_{in}^L|^2 1_{\{|Z_{in}^L| > \epsilon\}}] \le \epsilon^{-2} 
\E_{\frac{i-1}{n}}[|Z_{in}^L|^4].$$}
Then, by Markov inequality, it suffices 
to prove that 
{\rv
\begin{equation}
  \sum_{i = 1}^{\lfloor nt \rfloor} \E_{\frac{i-1}{n}}[|Z_{in}^L|^4] \xrightarrow{\mathbb{P}} 0 .
\label{eq: cond2 local}
\end{equation}
}
We have that 
\begin{align*}
\E_{\frac{i-1}{n}}[(Z_{in}^L)^4] \leq C\left( n^{2(\frac{1}{\alpha} - 1)} 
\E_{\frac{i-1}{n}}[f^4(n^{\frac{1}{\alpha}}(x - X_{\frac{i - 1}{n}}), n^{\frac{1}{\alpha}} \Delta_i^n X)] +n^{2(1 - \frac{1}{\alpha})}  
\E_{\frac{i-1}{n}}[L^4_{[\frac{i-1}{n}, \frac{i}{n}]}(x)] \right). 
\end{align*}
From Lemma \ref{l: moments L and O}(a) we conclude that 
\[
\sum_{i = 1}^{\lfloor nt \rfloor} \E_{\frac{i-1}{n}}[|Z_{in}^L|^4] 
\leq C n^{2(\frac{1}{\alpha} - 1)} \sum_{i = 1}^{\lfloor nt \rfloor} 
h\left( n^{\frac{1}{\alpha}} (X_{\frac{i-1}{n}} - x) \right)
\]
where $h(y)=\E[f^4(y,X_1)+L_{[0,1]}(y)^4]$. It is easy to check that  {\rv $h\in L^1(\R)$}. Indeed, similarly as in the proof of $\psi \in L^1(\R)$, Minkowski and Jensen inequalities imply 
$$h(y) \le C \E[L_{[0, 1]}^4(y)] \le C \E[L_{[0, 1]}^4(0)] \mathbb{P}[\tau_y < 1].$$
Then the boundedness of moments of the local time in $0$ together with \eqref{eq: prob passage time} provides $h \in L^1(\R)$. As $h^2(y) \le C \E[L_{[0, 1]}^8(y)]$, following the same route it is easy to see that $h^2 \in L^1(\R)$.
Since $\alpha>1$, we deduce by Theorem  \ref{thm:Jeganathan} that 
\[
n^{2(\frac{1}{\alpha} - 1)} \sum_{i = 1}^{\lfloor nt \rfloor} 
h\left( n^{\frac{1}{\alpha}} (X_{\frac{i-1}{n}} - x) \right) \xrightarrow{\mathbb{P}} 0.
\]
This concludes the proof of Theorem \ref{th: main local} for local times. \\ \\

\noindent Now we proceed to the analysis of the occupation time. 
%
As for the local time case, the proof is based on a martingale approach. The definition of $W_{O,t}^n$ together with the approximation of $S_{O,t}^n$ as in Proposition \ref{l: estim local} provides 
\begin{equation}
W_{O,t}^n = n^{\frac{1}{2}(1 + \frac{1}{\alpha})}(S_{O,t}^n - O_{[0, t]}(x))
= \sum_{i = 1}^{\lfloor nt \rfloor} Z_{in}^{O}  - n^{\frac{1}{2}(1 + \frac{1}{\alpha})} {N}_{O,t}^{ n} + \mathcal{E}_{O,t}^{n}(x),
\label{eq: 3.10.5}
\end{equation}
where $Z_{in}^{O}$ is the principal term, given by 
$$Z_{in}^{O}:= n^{\frac{1}{2}(\frac{1}{\alpha}-1)} F(n^{\frac{1}{\alpha}}(x - X_{\frac{i - 1}{n}}), n^{\frac{1}{\alpha}} \Delta_i^n X) - n^{\frac{1}{2}(1 + \frac{1}{\alpha})} O_{[\frac{i-1}{n}, \frac{i}{n}]}(x),$$
while 
$$N_{O,t}^{ n}:= O_{[\frac{\lfloor nt \rfloor}{n}, t]}(x).$$
In a similar way as for $N_{L,t}^{n}$ we can show that $N_{O,t}^{ n}$ is negligible. Indeed, the following lemma holds true. Its proof can be found in the next section {\rv and is a direct consequence of \eqref{eq: 1.4.5}}. 
\begin{lemma}
The process $N_O^{n}=\{N_{O,t}^{ n}\ ;\ t\geq0\}$ satisfies the following convergence  in probability uniformly over compact sets
$$\lim_n n^{\frac{1}{2}(1 + \frac{1}{\alpha})} N_O^{ n}  = 0. $$
\label{l: edge occ}
\end{lemma}

\noindent
Due to Proposition \ref{l: estim local} and Lemma \ref{l: edge occ} we can write  \eqref{eq: 3.10.5}  as 
$$W_{O,t}^n = \sum_{i = 1}^{\lfloor nt \rfloor} {Z}_{in}^{O} + o_{\mathbb{P}}(1).$$
We are dealing with martingale differences.  Indeed, $\E_{\frac{i-1}{n}}[{Z}_{in}^O] = 0$ because of Lemma \ref{l: moments L and O}(b) with $p=1$, the definition of $F$ and the independence of the increments of the process $X$. Therefore, similarly as before, our proof is based on \cite[Theorem 3-2]{Jacod}. In particular,  we want to show the following convergence statements:
\begin{equation}
\sum_{i = 1}^{\lfloor nt \rfloor} \E_{\frac{i-1}{n}}[|Z_{in}^{O}|^2] \xrightarrow{\mathbb{P}} {k}_{O}^2 L_{[0,t]}(x) \qquad \forall t \in [0, 1],
\label{eq: cond1 occ}
\end{equation}
\begin{equation}
\sum_{i = 1}^{\lfloor nt \rfloor} \E_{\frac{i-1}{n}}[Z_{in}^O \Delta_i^n M]  \xrightarrow{\mathbb{P}} 0 \qquad \forall t \in [0, 1], 
\label{eq: 3.15 Jacod occ}
\end{equation}
\begin{equation}
\sum_{i = 1}^{\lfloor nt \rfloor} \E_{\frac{i-1}{n}}[|Z_{in}^{O}|^2\Indi{\{ |Z_{in}^{O}| > \epsilon \}}] \xrightarrow{\mathbb{P}} 0 \qquad \forall \epsilon > 0.
\label{eq: cond2 occ}
\end{equation}
{\rv The condition gathered in \eqref{eq: 3.15 Jacod occ} should hold for all square integrable continuous martingales.} \\

\noindent We start by proving \eqref{eq: cond1 occ}. Similarly as for ${Z}_{in}^{O}$, from the definition of $F$ and Lemma \ref{l: moments L and O} (b) and (d) it follows
\begin{equation}{\label{eq: squared occ}}
\sum_{i = 1}^{\lfloor nt \rfloor} \E_{\frac{i-1}{n}}[ |{Z}_{in}^{O}|^2] = n^{\frac{1}{\alpha}-1}\sum_{i = 1}^{\lfloor nt \rfloor} \tilde{\psi}(n^{\frac{1}{\alpha}}(x - X_{\frac{i - 1}{n}})),
\end{equation}
with $\tilde{\psi}$ given by
\begin{align*}
\tilde{\psi}(z) &:= \E \big[(\E[O_{[0, 1]}(z) | X_1 ] - O_{[0,1]}(z))^2 \big].
\end{align*}
We want to show {\rv that $\tilde{\psi} \in L^1(\R)$}. Let $c_z$ be a deterministic constant only depending on $z$ and observe that by Jensen's inequality, 
it is enough to show the integrability of $\E[(O_{[0, 1]}(z)- c_z)^2]$, where $c_z$ is arbitrary. Take $c_z = 0$ for $z \ge 0$. Proceeding as in the case of local time, it is easy to see that 
$$\E[O^2_{[0, 1]}(z)] \le C \mathbb{P}(\tau_z < 1) \in L^1(\R).$$ 
For $z < 0$ we choose instead $c_z = 1$. Remarking that $\E[(1 - O_{[0, 1]}(z))^2] \le \mathbb{P}(\tau_z < 1)$ we obtain the same result. After that it is straightforward to prove that $\tilde{\psi}^2 \in L^1(\R)$, observing {\rv that} it is enough to show the integrability of  $\E[(O_{[0, 1]}(z)- c_z)^4]$. We can then apply Theorem \ref{thm:Jeganathan}, which implies \eqref{eq: cond1 occ} with $k_O^2 = \int_{\R} \tilde{\psi}(z) dz$. \\

 \noindent Condition \eqref{eq: 3.15 Jacod occ} is, as for the local time, a consequence of the martingale representation for jump measures in Lemma 3 (ii) and Theorem 6 of \cite{Repr}. \\

\noindent Finally, we show condition \eqref{eq: cond2 occ}. By  {\rv Cauchy-Schwarz} and Markov inequalities it is enough to prove that
$${\rv \sum_{i = 1}^{\lfloor nt \rfloor} \E_{\frac{i-1}{n}}[|{Z}_{in}^O|^4] \xrightarrow{\mathbb{P}} 0.}$$
By the definition of ${Z}_{in}^O$ and Lemma \ref{l: moments L and O},  following the same proof as to obtain \eqref{eq: squared occ}, it holds that
\begin{align*}
\sum_{i = 1}^{\lfloor nt \rfloor} \E_{\frac{i-1}{n}}[|{Z}_{in}^O|^4] = n^{2(\frac{1}{\alpha} - 1)} \sum_{i = 1}^{\lfloor nt \rfloor} \tilde{h}(n^{\frac{1}{\alpha}}(x - X_{\frac{i - 1}{n}})), 
\end{align*}
where $\tilde{h}(z) := \E \big[(\E[O_{[0, 1]}(z) | X_1 ] - O_{[0,1]}(z))^4 \big]$. Regarding the integrability of $\tilde{h}$ and $\tilde{h}^2$, proceeding as for $\tilde{\psi}$ it is sufficient to show the result for $\E \big[( O_{[0,1]}(z)- c_z)^4 \big]$ and $\E \big[( O_{[0,1]}(z)- c_z)^8 \big]$ respectively, with $c_z$ arbitrary. As before, the result follows by choosing $c_z =0$ for $z \ge 0$ and $c_z = 1$ for $z< 0$, and using \eqref{eq: prob passage time}. We deduce by Theorem \ref{thm:Jeganathan} that
\[
n^{2(\frac{1}{\alpha} - 1)} \sum_{i = 1}^{\lfloor nt \rfloor} 
\tilde{h}\left( n^{\frac{1}{\alpha}} (X_{\frac{i-1}{n}} - x) \right) \xrightarrow{\mathbb{P}} 0,
\]
which concludes the proof of our main theorem.
\end{proof}

\section{Proof of the auxiliary results}\label{sec4}
This section is devoted to the proof of technical lemmas and proposition we have previously stated and used in order to obtain our main result, gathered in Theorem \ref{th: main local}.  Before we  proceed we introduce the function $r(y) := e^{- \textbf{i} y}$, which will be useful in the sequel.

{\rv 
\subsection{Proof of Lemma \ref{l:lemmaFourierLT}}
\begin{proof}
The Fourier representation is proven in \cite[Proposition 11]{J04}, we here prove the boundedness of the moments of $L_t(x)$. From the Fourier representation stated in the first part of this lemma we have, for any $p \in \mathbb{N}$, 
\begin{multline}{\label{eq: local bounded 1}}
\E[(L_{[0, 1]}( x))^p]\\
\begin{aligned}
  & = \E\left[\left(\frac{1}{2 \pi} \int_{\R} \int_{[0,1]} r({\xi}_1 (X_{u_1}- x) du_1 d \xi_1\right) \times .. \times \left(\frac{1}{2 \pi} \int_{\R} \int_{[0,1]} r({\xi}_p (X_{u_p}- x) du_p d\xi_p\right)\right] \\
&= \frac{1}{(2 \pi)^p} \int_{\R^p} \int_{[0,1]^p} r(-({\xi}_1 +... + {\xi}_p)x) \E[r({\xi}_1 X_{u_1}) \times ... \times r({\xi}_p X_{u_p})] d\bm{u} d\bm{{\xi}}. 
\end{aligned}
\end{multline}
The Fubini theorem, applied in the last identity, will be justified once we show that $|\E[r({\xi}_1 X_{u_1}) \times ... \times r({\xi}_p X_{u_p})]|$ is integrable, which will be proved next. 
We can now assume without losing of generality that $p$ is even and $u_1 \le u_2 \le ... \le u_p$. Then, the expectation above can be seen as{
\begin{align*}
&\E[r((\xi_1 + ... + \xi_p) X_{u_1}) r((\xi_2 + ... + \xi_p) (X_{u_2} -X_{u_1})) \times ... \times r(\xi_p(X_{u_p} -X_{u_{p - 1}}))] \\
& = \E[r((\xi_1 + ... + \xi_p) X_{u_1})] \times ... \times \E[r(\xi_p(X_{u_p} -X_{u_{p - 1}}))]
\end{align*}
having employed the self-similarity of $X$ and the stationarity and independence of its increments. Relation \eqref{eq:characteristifcuntion} then implies that
\begin{align*}
&|\E[r((\xi_1 + ... + \xi_p) X_{u_1}) r((\xi_2 + ... + \xi_p) (X_{u_2} -X_{u_1})) \times ... \times r(\xi_p(X_{u_p} -X_{u_{p - 1}}))]|\\
&=e^{- |\xi_1 + ... + \xi_p|^\alpha u_1} \times ... \times e^{- |\xi_p|^\alpha (u_p - u_{p-1})},
\end{align*}
}
Now, up to apply the change of variables $\xi_j + ... + \xi_p =: \eta_j$ for any $j =1, ... , p$, it is easy to check that the {absolute value in the integrand of \eqref{eq: local bounded 1} is integrable, thus yielding the finiteness of the $p$-moment of $L_{t}$.}

\end{proof}
}
\subsection{Proof of Proposition \ref{l: estim local}}
\begin{proof}
Recall that $\mathcal{A}_n:=\sigma(\{X_{\frac{i}{n}}\ ;\ i\in\N\})\subset\mathcal{F}$. We first deal with the analysis of the local time. Observe that the problem is reduced to proving the following two claims:
\begin{enumerate}
\item[(i)] First of all show that 
\begin{equation}
\E\left[L_{[\frac{i-1}{n}, \frac{i}{n}]}(x) | X_{\frac{i-1}{n}},  \Delta_i^n X\right] = n^{\frac{1}{\alpha} - 1} f\left(n^{\frac{1}{\alpha}}(x - X_{\frac{i - 1}{n}}), n^{\frac{1}{\alpha}} \Delta_i^n X\right),
\label{eq: local intervals}
\end{equation}
for {\rv $f(x, y)= \E[L_{[0,1]}(x) | X_1 = y]$.}
\item[(ii)] After that we prove that $\mathcal{E}_L^{(n)}(x)$ satisfies
\begin{align}\label{eq: edge local} 
\sup_{\rv 0\leq t \leq T}n^{\frac{1}{2}(1-\frac{1}{\alpha})}|\mathcal{E}_{L,t}^{(n)}(x)| \xrightarrow{\mathbb{P}}0,
\end{align}
for all $T>0.$
\end{enumerate}
If (i) and (ii) are satisfied, it is easy to check that the result on the local time stated in Proposition \ref{l: estim local} holds true. To prove this reduction we observe that by the independent increments property of $X$, 
\begin{align*}
S_{L,t}^{(n)}(x)
  &=\E\left[L_{[ \lfloor nt\rfloor/n,t]}{\rv (x)}\ |\ \mathcal{A}_{n}\right]+\sum_{k=1}^{\lfloor nt\rfloor}\E\left[L_{[\frac{k-1}{n},\frac{k}{n}]}{\rv (x)}\ |\ \mathcal{A}_{n}\right]\\
	&=\E\left[L_{[ \lfloor nt\rfloor/n,t]}{\rv (x)}\ |\ 
 {\rv \mathcal{A}_{n} }\right]+\sum_{k=1}^{\lfloor nt\rfloor}\E\left[L_{[\frac{k-1}{n},\frac{k}{n}]}{\rv (x)}\ |\  X_{\frac{k-1}{n}},\Delta_{k}^nX\right],
\end{align*}
so that
relation {\rv \eqref{eq: local intervals}} implies the desired result. We now proceed with the proof of \eqref{eq: local intervals} and \eqref{eq: edge local}.  In order to show \eqref{eq: local intervals}, {we consider the approximated local time $L_{t}^{\varepsilon}(x),$
defined by 
$$L_{t}^{\varepsilon}(x)
  :=\int_0^t\phi_{\varepsilon}(X_{s}-x)ds,$$ 
where $\phi_{\varepsilon}(x):=(2\pi\varepsilon)^{-1/2}\exp\{-\frac{1}{2\varepsilon}x^2\}$. By \cite{J04}, we have the convergence towards zero of $\|L_{t}(x)-L_{t}^{\varepsilon}(x)\|_{L^{2}(\Omega)}\rightarrow0$, which   implies that 
\begin{align*}
\E[L_{[\frac{i-1}{n},\frac{i}{n}]}(x) \ |\ X_{\frac{i-1}{n}},\Delta_{i}^nX]
  &=\lim_{\varepsilon\rightarrow0}\E[L^{\varepsilon}_{[\frac{i-1}{n},\frac{i}{n}]} (x) \ |\ X_{\frac{i-1}{n}},\Delta_{i}^nX],
\end{align*}
where $L_{[a,b]}^{\varepsilon}(x):=L_{b}^{\varepsilon}(x)-L_{a}^{\varepsilon}(x)$. By the self-similarity of $X$, for every $a,b\in\R$,
\begin{multline*}
\E[L^{\varepsilon}_{[\frac{i-1}{n},\frac{i}{n}]}(x)\ |\ X_{\frac{i-1}{n}}=a, \Delta_{i}^nX=b]\\
\begin{aligned}
  &=
\E[\int_{\frac{i-1}{n}}^{\frac{i}{n}}\phi_{\varepsilon}(X_{s}-x)ds\ |\ X_{\frac{i-1}{n}}=a, \Delta_{i}^nX=b]\\
&=
\E[\int_{\frac{i-1}{n}}^{\frac{i}{n}}\phi_{\varepsilon}(n^{-1/\alpha}X_{ns}-x)ds\ |\ n^{-1/\alpha}X_{ i-1}=a, n^{-1/\alpha}(X_{i}-X_{i-1})=b]\\
&=
\int_{\frac{i-1}{n}}^{\frac{i}{n}}\E[\phi_{\varepsilon}(n^{-1/\alpha}(X_{ns}-X_{i-1})+a-x)\ |\ n^{-1/\alpha}X_{i-1}=a, n^{-1/\alpha}(X_{i}-X_{i-1})=b]ds,
\end{aligned}
\end{multline*}
where the last identity follows by Fubini's theorem, which is applicable due to the fact that $\phi_{\varepsilon}$ is bounded and the domain of integration is compact. The independent increments property of $X$ implies that $X_{ i-1}$ is independent of the vector $(X_{s}-X_{i-1},X_{i}-X_{i-1})$. Moreover, by the stationarity of the increments of $X$, we have that $$(X_{s}-X_{i-1},X_{i}-X_{i-1})\stackrel{Law}{=}(X_{s-i+1},X_{1}).$$ 
This identity combined with a suitable change of variables yields
\begin{multline*}
\E[L^{\varepsilon}_{[\frac{i-1}{n},\frac{i}{n}]}(x)\ |\ X_{\frac{i-1}{n}}=a,n^{1/\alpha}\Delta_{i}^nX=b]\\
\begin{aligned}
  &=
n^{-1}\int_{0}^{1}\E[\phi_{\varepsilon}(n^{-1/\alpha}X_{u}+a-x)\ |\ X_1=n^{1/\alpha}b]ds\\
  &=
n^{1/\alpha-1}\int_{0}^{1}\E[\phi_{\varepsilon n^{-2/\alpha}}( X_{u}+an^{1/\alpha}-xn^{1/\alpha})\ |\ X_1=n^{1/\alpha}b]ds\\
  &=
n^{1/\alpha-1}\E[\int_{0}^{1}\phi_{\varepsilon n^{-2/\alpha}}(X_{u}+an^{1/\alpha}-xn^{1/\alpha})ds\ |\ X_1=n^{1/\alpha}b].
\end{aligned}
\end{multline*}
Using the fact that 
$$\int_{0}^{1}\phi_{\varepsilon n^{-2/\alpha}}( X_{u}+an^{1/\alpha}-xn^{1/\alpha})du
\stackrel{L^{2}(\Omega)}{\rightarrow}L_{1}(n^{1/\alpha}(x-a)),$$
as $\varepsilon\rightarrow0$, we obtain 
\begin{align*}
\lim_{\varepsilon\rightarrow0}\E[L^{\varepsilon}_{[\frac{i-1}{n},\frac{i}{n}]}(x) \ |\ X_{\frac{i-1}{n}}=a,n^{1/\alpha}\Delta_{i}^nX=b]
  &= n^{1/\alpha-1}\E[L_{1}(n^{1/\alpha}(x-a))\ |\ X_{1}=n^{1/\alpha}b].
\end{align*}
This finishes the proof of \eqref{eq: local intervals}.\\
}\\

\noindent We are left with the problem of showing \eqref{eq: edge local}. 
\noindent Due to \eqref{eq: 1.4.5}, for every $T,\varepsilon> 0$ and $x \in \R$, we have
$$\sup_{0\leq t\leq T}|L_{[\frac{\lfloor nt \rfloor}{n}, t]}(x)| \le C n^{-1 + \frac{1}{\alpha} + \epsilon},$$
for any $\epsilon \in (0, 1 - 1/\alpha)$. 
Relation \eqref{eq: edge local} follows from here. \\

\noindent Next we deal with identity \eqref{eq:SOexact}. By \eqref{eq:Oabtolocaltime} and \eqref{eq: local intervals}, we have that 
\begin{align*}
\E[O_{t}(x)\ |\ \mathcal{A}_{n}] &=\sum_{i=1}^{\lfloor nt\rfloor} \E[O_{[\frac{i-1}{n}, \frac{i}{n}]}(x)\ |\ \mathcal{A}_{n}] + \E\left[O_{[\frac{\lfloor nt \rfloor}{n},t]}(x)\ |\ \mathcal{A}_n\right]\\
  &=\sum_{i=1}^{\lfloor nt\rfloor} \int_{x}^{\infty}\E[L_{[\frac{i-1}{n}, \frac{i}{n}]}(y)\ |\ \mathcal{A}_{n}]dy + \E\left[O_{[\frac{\lfloor nt \rfloor}{n},t]}(x)\ |\ \mathcal{A}_n\right]\\
  &=n^{\frac{1}{\alpha}-1}\int_{x}^{\infty} \sum_{i=1}^{\lfloor nt\rfloor}f\left(n^{\frac{1}{\alpha}}(y-X_{\frac{i-1}{n}}),n^{\frac{1}{\alpha}}\Delta_{i}^{n}X\right)dy+\E\left[O_{[\frac{\lfloor nt \rfloor}{n},t]}(x)\ |\ \mathcal{A}_n\right]\\
  &=\frac{1}{n} \sum_{i=1}^{\lfloor nt\rfloor} F\left(n^{\frac{1}{\alpha}}(x-X_{\frac{i-1}{n}}),n^{\frac{1}{\alpha}}\Delta_{i}^{n}X\right)+ \E\left[O_{[\frac{\lfloor nt \rfloor}{n},t]}(x)\ |\ \mathcal{A}_n\right],
\end{align*}
where the last identity follows from a suitable change of variables. {Consequently, we deduce \eqref{eq:SOexact}}. In addition, {\rv by definition of occupation $O_{[a,b]}(x)$} we have that 
$|\E[O_{[\frac{\lfloor nt \rfloor}{n},t]}(x)\ |\ \mathcal{A}_n]|\leq n^{-1}$ which implies that
$$\sup_{0\leq t\leq T}n^{\frac{1}{2}(1+\frac{1}{\alpha})}\E\left[O_{[\frac{\lfloor nt \rfloor}{n},t]}(x)\ |\ \mathcal{A}_n\right]\leq n^{\frac{1}{2}(1/\alpha-1)},$$ 
which converges towards zero. 
{\rv The proof is now complete.}
\end{proof}

\subsection{Proof of Lemma \ref{l: moments L and O}}
\begin{proof}
\textit{Part (a)} \\
Equation \eqref{eq: local intervals} provides the desired result for $p = 1$  as the increments are independent and so, in particular, $X_{\frac{i-1}{n}}$ is independent from $\Delta_i^n X$. We thus study the case $p \ge 2$. 
The representation of the local time as in Lemma \ref{l:lemmaFourierLT} leads us to 
\begin{align*}
 \E_{\frac{i-1}{n}}\left[L^p_{[\frac{i-1}{n}, \frac{i}{n}]}(x)\right]
= & \E_{\frac{i-1}{n}}\left[\frac{1}{(2 \pi)^p} \int_{\R^p} \int_{[\frac{i-1}{n}, \frac{i}{n}]^p} r(\xi_1(X_{s_1} - x)) ... r(\xi_p(X_{s_p} - x)) d\bm{s} d\bm{\xi}\right],
\end{align*}
where we recall that $\E_{\frac{i-1}{n}}[\cdot] = \E[\cdot|\mathcal{F}_{\frac{i-1}{n}}]$ and $r(y) = e^{- \bm{i} y}$. \\
The change of variables $s_j =: \frac{i-1}{n} + \frac{u_j}{n}$ for $j=1,... , p$ yields
\begin{align*}
 & \frac{1}{n^p}  \frac{1}{(2 \pi)^p} \int_{\R^p} \int_{[0,1]^p} r((\xi_1 + ... + \xi_p)(X_{\frac{i-1}{n}} - x)) \E_{\frac{i-1}{n}}[r(\xi_1(X_{\frac{i-1}{n} + \frac{u_1}{n}} - X_{\frac{i-1}{n}})) \times ...\\
& \times r(\xi_p(X_{\frac{i-1}{n} + \frac{u_p}{n}} - X_{\frac{i-1}{n}}))] d\bm{u} d\bm{\xi} \\
= & n^{-p} \frac{1}{(2 \pi)^p} \int_{\R^p} \int_{[0,1]^p} r((\xi_1 +... +  \xi_p)(X_{\frac{i-1}{n}} - x)) \E[r(\xi_1 X_{\frac{u_1}{n}}) \times ... \times r(\xi_p X_{ \frac{u_p}{n}})] d\bm{u} d\bm{\xi},
\end{align*}
having also used the independence {and stationarity of the increments of $X$, as well as Fubini's theorem, which is justified as in the proof of Lemma \ref{l:lemmaFourierLT}}. As the process $X$ is self-similar, we obtain 
\begin{align*}
n^{-p } \frac{1}{(2 \pi)^p} \int_{\R^p} \int_{[0,1]^p} r((\xi_1 +... + \xi_p)(X_{\frac{i-1}{n}} - x)) \E[r(\xi_1 n^{-\frac{1}{\alpha}} X_{u_1})\times ... \times r(\xi_p n^{-\frac{1}{\alpha}} X_{u_p})] d\bm{u} d\bm{\xi}.
\end{align*}
Applying the change of variables $\tilde{\xi}_j = n^{-\frac{1}{\alpha}} \xi_j$ for $j=1,... , p$, we get 
\begin{align*}
& n^{p( \frac{1}{\alpha}- 1)} \frac{1}{(2 \pi)^p} \int_{\R^p} \int_{[0,1]^p} r((\tilde{\xi}_1 + ... + \tilde{\xi}_p)n^{\frac{1}{\alpha}}(X_{\frac{i-1}{n}} - x)) \E[r(\tilde{\xi}_1 X_{u_1}) \times ... \times r(\tilde{\xi}_p X_{u_p})] d\bm{u} d\bm{\tilde{\xi}} \\
  &=  n^{p(\frac{1}{\alpha}-1)}  \varphi_1(n^{\frac{1}{\alpha}}(x - X_{\frac{i - 1}{n}})).
\end{align*}
{\rv It is indeed 
\begin{align*}
 &\frac{1}{(2 \pi)^p} \int_{\R^p} \int_{[0,1]^p} r(-({\xi}_1 +... + {\xi}_p)y) \E[r({\xi}_1 X_{u_1}) \times ... \times r({\xi}_p X_{u_p})] d\bm{u} d\bm{{\xi}} \\
& =\E[(L_{[0, 1]}( y))^p] = \varphi_1(y),
\end{align*}
because of \eqref{eq: local bounded 1}.}
It is important to remark that, in the proof above, it is possible to obtain the {\rv desired}
result as we consider the conditional expectation with respect to $\mathcal{F}_{\frac{i-1}{n}}$, and we used several times the independence of the increments. \\

\noindent \textit{Part (b)}\\
The definition of occupation time provides
\begin{align*}
& \E_{\frac{i-1}{n}}\left[O^p_{[\frac{i-1}{n}, \frac{i}{n}]}(x)\right] \\
& = \E_{\frac{i-1}{n}}\left[\left(\int_{x}^\infty L_{[\frac{i-1}{n}, \frac{i}{n}]}(y_1) dy_1\right) \times ... \times \left(\int_{x}^\infty L_{[\frac{i-1}{n}, \frac{i}{n}]}(y_p) dy_p\right) \right] \\
& =  \int_{x}^\infty ... \int_{x}^\infty \E_{\frac{i-1}{n}}\left[ L_{[\frac{i-1}{n}, \frac{i}{n}]}(y_1) \times ... \times  L_{[\frac{i-1}{n}, \frac{i}{n}]}(y_p)\right] d\bm{y}. 
\end{align*}
{The change in the order of integration is justified by Tonelli's theorem, which is valid due to the fact that $L_{[a,b]}(x)\geq0$ for all $a\leq b$ and $x\geq0$.\\
Acting as in the proof of Part (a) it is then easy to check that
\begin{align*}
&\E_{\frac{i-1}{n}}\left[O^p_{[\frac{i-1}{n}, \frac{i}{n}]}(x)\right] \\
= & \frac{n^{p(\frac{1}{\alpha} - 1)}}{(2 \pi)^p} \int_{[x, \infty]^p} \int_{\R^p} \int_{[0, 1]^p} r(\tilde{\xi}_1 n^{\frac{1}{\alpha}}(X_{\frac{i-1}{n}}- y_1)) \times ... \times r(\tilde{\xi}_p n^{\frac{1}{\alpha}}(X_{\frac{i-1}{n}}- y_p)) \\
& \times \E\left[r(\tilde{\xi}_1 X_{u_1}) \times ... \times r(\tilde{\xi}_p X_{u_p})\right] d\bm{u} d\bm{\tilde{\xi}} d\bm{y}.
\end{align*}}
To conclude the analysis we apply the change of variable $\tilde{y}_j := n^{\frac{1}{\alpha}}(y_j -X_{\frac{i-1}{n}})$ for $j=1,..., p$. We get
\begin{align*}
& \frac{n^{-p}}{(2 \pi)^p} \int_{[n^{\frac{1}{\alpha}}(x - X_{\frac{i-1}{n}}), \infty]^p} \int_{\R^p} \int_{[0, 1]^p} r(-\tilde{\xi}_1 \tilde{y}_1) \times ... \times r(-\tilde{\xi}_p \tilde{y}_p) \\
& \times \E[r(\tilde{\xi}_1 X_{u_1}) \times ... \times r(\tilde{\xi}_p X_{u_p})] d\bm{u} d\bm{\tilde{\xi}} d\bm{\tilde{y}} \\
=& n^{ - p} \varphi_2(n^{\frac{1}{\alpha}}(x - X_{\frac{i-1}{n}})),
\end{align*}
with 
\begin{align*}
\varphi_2(z) &= \int_{[z, \infty]^p}  \frac{1}{(2 \pi)^p} \int_{\R^p} \int_{[0, 1]^p} \E[r({\xi}_1 (X_{u_1}- {y}_1)) \times ... \times r({\xi}_p (X_{u_p}- {y}_p))] d\bm{u} d\bm{{\xi}} d\bm{y} \\
& = \int_{[z, \infty]^p} \E[L_{[0, 1]}(y_1)\times ... \times L_{[0, 1]}(y_p)] dy_1 ... dy_p 
= \E[O_{[0,1]}^p(z)],
\end{align*}
{where the previous to last identity is justified by Tonelli's theorem.} The proof of this part is therefore completed. \\

\noindent \textit{Part (c)}\\
{\rv According to} \eqref{eq: local intervals} it is 
$$ f\left(n^{\frac{1}{\alpha}}(x - X_{\frac{i - 1}{n}}), n^{\frac{1}{\alpha}} \Delta_i^n X\right) = n^{1 - \frac{1}{\alpha}} \E\left[L_{[\frac{i-1}{n}, \frac{i}{n}]}(x) | X_{\frac{i-1}{n}}, \Delta_i^n X\right].$$
Then, following the proof of Part (a), one can easily {\rv obtain}
\begin{equation}{\label{eq: A}}
f\left(n^{\frac{1}{\alpha}}(x - X_{\frac{i - 1}{n}}), n^{\frac{1}{\alpha}} \Delta_i^n X\right) = \E\left[L_{[0,1]}(n^{\frac{1}{\alpha}}(x - X_{\frac{i - 1}{n}})) | X_1\right].
\end{equation}
From the representation of the local time as in Lemma \ref{l:lemmaFourierLT} {\rv we conclude the
identity}
\begin{align*}
&\E_{\frac{i - 1}{n}} \left[ f\left(n^{\frac{1}{\alpha}}(x - X_{\frac{i - 1}{n}}), n^{\frac{1}{\alpha}} \Delta_i^n X\right) L_{[\frac{i-1}{n}, \frac{i}{n}]}(x)  \right] \\
& = \frac{1}{2 \pi} \int_{\R} \int_{\frac{i - 1}{n}}^{\frac{i}{n}} \E_{\frac{i - 1}{n}} \left[ \E\left[L_{[0,1]}(n^{\frac{1}{\alpha}}(x - X_{\frac{i - 1}{n}})) | X_1\right] r(\xi(X_s -x))  \right] ds d\xi.
\end{align*}
The change of variable $s:= \frac{i-1}{n} + \frac{u}{n}$ provides the quantity above is equal to
\begin{align*}
& \frac{1}{n} \frac{1}{2 \pi} \int_{\R} \int_{0}^{1} r(\xi(X_{\frac{i-1}{n}} -x)) \E_{\frac{i - 1}{n}} \left[ \E\left[L_{[0,1]}(n^{\frac{1}{\alpha}}(x - X_{\frac{i - 1}{n}})) | X_1\right] r(\xi(X_{\frac{i-1}{n} + \frac{u}{n} } -X_{\frac{i-1}{n}}))  \right] du d\xi \\
& = n^{\frac{1}{\alpha}-1} \frac{1}{2 \pi} \int_{\R} \int_{0}^{1} r({\rv \tilde{\xi}} n^\frac{1}{\alpha}(X_{\frac{i-1}{n}} -x)) \E \left[ \E\left[L_{[0,1]}(n^{\frac{1}{\alpha}}(x - X_{\frac{i - 1}{n}})) | X_1\right] r(\tilde{\xi}X_{u})  \right] du d\tilde{\xi},
\end{align*}
having used the independence of the increments, the self-similarity of the process $X$ and the change of variable $\tilde{\xi} := n^{-\frac{1}{\alpha}} \xi$. The proof is concluded once we observe that this is $n^{\frac{1}{\alpha}- 1} \varphi_3(n^\frac{1}{\alpha}(x - X_{\frac{i-1}{n}}))$, {\rv as
\begin{align*}
 & \frac{1}{2 \pi} \int_{\R} \int_{0}^{1} \E \left[ \E\left[L_{[0,1]}(z) | X_1\right] r(\tilde{\xi}(X_{u}-z))  \right] du d\tilde{\xi} \\
& = \E \left[ \E\left[L_{[0,1]}(z) | X_1\right] L_{[0,1]}(z)  \right] = \varphi_3(z).
\end{align*}}
\vspace{0.3 cm}

\noindent \textit{Part (d)}\\
We remark that, from the definition of $F$, Equation \eqref{eq: A} and a suitable change of variable, it is 
\begin{align*}
F\left(n^{\frac{1}{\alpha}}(x - X_{\frac{i - 1}{n}}), n^{\frac{1}{\alpha}} \Delta_i^n X\right) & = n^{\frac{1}{\alpha}} \int_x^\infty f\left(n^{\frac{1}{\alpha}}(y - X_{\frac{i - 1}{n}}), n^{\frac{1}{\alpha}} \Delta_i^n X\right)dy \\
& = n^{\frac{1}{\alpha}} \int_x^\infty \E\left[L_{[0,1]}(n^{\frac{1}{\alpha}}(y - X_{\frac{i - 1}{n}})) | X_1\right] dy \\
& = \E\left[O_{[0,1]}(n^{\frac{1}{\alpha}}(x - X_{\frac{i - 1}{n}})) | X_1\right].
\end{align*}
Then, following the same route as in the proof of {\rv parts} (b) and (c), it is easy to check that
\begin{align*}
&\E \left[F\left(n^{\frac{1}{\alpha}}(x - X_{\frac{i - 1}{n}}), n^{\frac{1}{\alpha}} \Delta_i^n X \right)
    O_{[\frac{i-1}{n}, \frac{i}{n}]}(x) | \mathcal{F}_{\frac{i - 1}{n}} \right]\\
    &= n^{ - 1} \E \left[O_{[0, 1]} (n^{\frac{1}{\alpha}}(x - X_{\frac{i-1}{n}})) \E[O_{[0, 1]} (n^{\frac{1}{\alpha}}(x - X_{\frac{i-1}{n}}))| X_1] \right] \\
    & {\rv = n^{ - 1} \varphi_4(n^{\frac{1}{\alpha}}(x - X_{\frac{i-1}{n}}),}
\end{align*}
as {\rv required}.

\end{proof}

\subsection{Proof of Lemmas \ref{l: edge local} and \ref{l: edge occ}}
\begin{proof}
It follows by arguments analogous to part (ii) in Proposition \ref{l: estim local}. In particular, the statement of Lemma \ref{l: edge local} is a straightforward consequence of the {\rv H\"older} property in \eqref{eq: 1.4.5}.

\end{proof}

\bibliographystyle{chicago}
 
\end{document}